\documentclass[11pt]{amsart}
\usepackage[b5paper, left=0.855in, right=0.855in, top=1.1in, bottom=1.35in]{geometry}
\usepackage{xcolor}
\usepackage{colortbl}

\usepackage{xpatch}
\usepackage{tabu}

\makeatletter
\xpatchcmd{%
\@maketitle}{%
\ifx\@empty\authors \else \@setauthors \fi
}{%
  \ifx\@empty\authors \else \@setauthors \fi

  \ifx\@empty\addresses \else\@setaddresses\fi
}{\typeout{Patch successful}}{\typeout{Patch failed}}
\makeatother

\usepackage{amsmath}
\usepackage{amssymb}
\usepackage{mathtools}
\usepackage{hhline}

\usepackage[bookmarks=false]{hyperref}

\newcommand{\var}{\mathit{var}}
\newcommand{\true}{1}
\newcommand{\false}{0}
\newcommand{\drat}{\mathrm{DRAT}}
\title{Computing Small Unit-Distance Graphs \mbox{with Chromatic Number 5}}
\author{Marijn J.H. Heule}
\address{Department of Computer Science, The University of Texas at Austin}
\email{marijn@heule.nl}

\thanks{
 Supported by the National Science Foundation (NSF) under grant CCF-1526760.}
 
\begin{document}
\maketitle
\begin{abstract}
We present a new method for reducing the size of graphs with a given property.
Our method, which is based on clausal proof minimization, allowed us to compute
several 553-vertex unit-distance graphs with chromatic number 5, while the
smallest published unit-distance graph with chromatic number 5 has 1581
vertices.
The latter graph
was constructed by Aubrey de Grey to show that the chromatic number of the plane is at least 5.
The lack of a 4-coloring of our graphs is due to a clear pattern enforced on some vertices. 
Also, our graphs can be mechanically validated in a second, which suggests that the pattern is based on a
reasonably short argument. 
\end{abstract}

\section{Introduction}
The {\em chromatic number of the plane}, a problem first proposed by 
Edward Nelson in 1950~\cite{coloring},
asks how many colors are needed to color all points of the plane such that no two points at distance 1 from each other have the same color.
Early results showed that at least four and at most seven colors are required.
By the de Bruijn--Erd\H{o}s theorem, the chromatic number of the plane is the largest possible chromatic number
of a finite unit-distance graph~\cite{deBruijn}.
The Moser Spindle, a unit-distance graph with 7 vertices and 11 edges, shows the lower bound~\cite{Moser},
while the upper bound is due to a 7-coloring of the entire plane by John Isbell~\cite{coloring}.

In a recent breakthrough for this problem, Aubrey de Grey improved the lower bound by providing
a unit-distance graph with 1581 vertices with chromatic number 5~\cite{DeGrey}. This graph was obtained 
by shrinking the initial graph with chromatic number 5 consisting of $20\,425$ vertices. The 1581-vertex graph is 
almost minimal: at most 4 vertices can be removed without introducing a 4-coloring of the remaining graph. 
The discovery by de Grey started a Polymath project to find smaller unit-distance graphs with chromatic number 5.

In this paper, we present our first contributions in this direction and we describe a method
to reduce the size of graphs while preserving their chromatic number. Our results show that
this method is quite effective  as it was able to produce unit-distance graphs with 553 vertices. Our method 
exploits two formal-methods technologies: the ability of  \emph{satisfiability \textup{(}SAT\textup{)}
solvers} to find a short refutation for unsatisfiable formulas (if they exist) and \emph{proof checkers} that can minimize refutations and unsatisfiable formulas. 

The refutations emitted by SAT solvers are hardly minimal. Depending on the application from which the 
formula originates, typically $10\%$ to $99\%$ of the refutation can be omitted. Several techniques have 
been developed to avoid checking irrelevant parts of a refutation~\cite{Heule:2013:trim}. These techniques
minimize proofs in order to share and revalidate them.
For example, the proofs of the Boolean Pythagorean Triples~\cite{ptn} and Schur Number Five~\cite{S5}
problems are enormous, even after minimization: 200 terabytes and 2 petabytes, respectively. 

Here we use clausal-proof-minimization techniques for a different purpose: shrinking graphs.
Given a unit-distance graph with chromatic number 5, we first construct a propositional formula that
encodes whether there exists a valid 4-coloring of this graph. This formula is unsatisfiable and 
we can use a SAT solver to compute a refutation. From the minimized
refutation, we extract a subgraph that also has chromatic number 5. 
We then apply this process repeatedly to make the graph ever smaller.

\section{Preliminaries}

\subsection{Chromatic Number of the Plane}

The Chromatic Number of the Plane (CNP) asks how many colors are required in a coloring of the plane to
ensure that there exists no monochromatic pair of points with distance 1. A graph for which all edges
have the same length is called a {\em unit-distance graph}. A lower bound for CNP of $k$
colors can be obtained by showing that a unit-distance graph has chromatic number $k$.

We will use three operations to construct larger and larger graphs:
the Minkowski sum, rotation, and merge.
Given two sets of points $A$ and $B$, the Minkowski sum of $A$ and $B$, denoted by $A \oplus B$, equals
\mbox{$\{a+b \mid a \in A, b \in B\}$}. Consider the sets of points $A = \{(0,0), (1,0)\}$ and \mbox{$B = \{(0,0), (1/2,\sqrt{3}/2)\}$}, then 
$A \oplus B = \{(0,0), (1,0), (1/2,\sqrt{3}/2), (3/2,\sqrt{3}/2)\}$.

\begin{figure}[t]
\includegraphics[width=\textwidth]{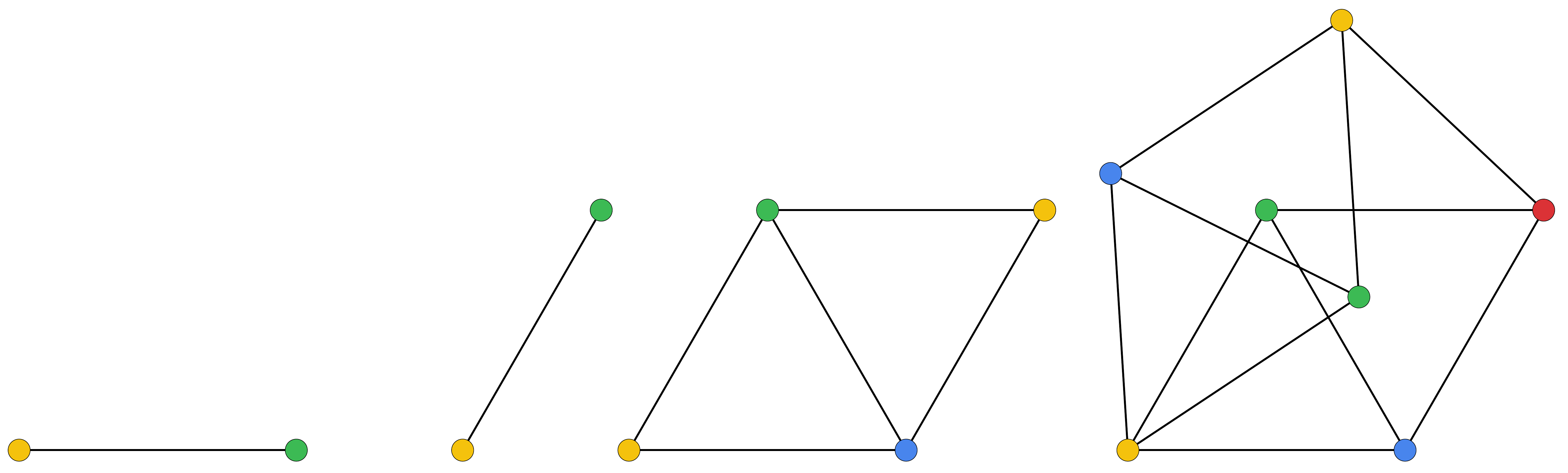}
\caption{From left to right: illustrations of $A$, $B$, $A \oplus B$, and the Moser Spindle.
The graphs shown have chromatic number 2, 2, 3, and 4, respectively.}
\label{fig:intro}
\end{figure}

Given a positive integer $i$, we denote by $\theta_i$ the rotation around point $(0,0)$ with angle $\arccos(\frac{2i-1}{2i})$
and by $\theta_i^k$ the application of $\theta_i$ $k$ times.
Let $p$ be a point with distance $\sqrt{i}$ from $(0,0)$, then the
points $p$ and $\theta_i(p)$ are exactly distance 1 apart and thus would be connected with an 
edge in a unit-distance graph. Consider again the set of points
$A \oplus B$ above. 
The points $A \oplus B \cup \theta_3(A \oplus B)$ form the Moser Spindle. 
Figure~\ref{fig:intro} shows visualizations of these sets with connected vertices colored differently. 

\subsection{Propositional Formulas}

We will minimize graphs on the propositional level. 
We consider propositional formulas in \emph{conjunctive normal form} (CNF), 
which are defined as follows. 
A \emph{literal} is either a variable $x$ (a \emph{positive literal}) 
or the negation $\overline x$ of a variable~$x$ (a \emph{negative literal}). 
The \emph{complement} $\overline l$ of a literal $l$ is defined as 
$\overline l = \overline x$ if $l = x$ and $\overline l = x$ if $l = \overline x$.
For a literal $l$, $\var(l)$ denotes the variable of $l$.
A \emph{clause} is a disjunction of literals and a \emph{formula} is a conjunction of clauses.

An \emph{assignment} is a function from a set of variables to the truth values 
\true{}~(\emph{true}) and \false{} (\emph{false}).
A literal $l$ is \emph{satisfied} by an assignment $\alpha$ if 
$l$ is positive and \mbox{$\alpha(\var(l)) = \true$} or if it is negative and $\alpha(\var(l)) = \false$.
A literal is \emph{falsified} by an assignment if its complement is satisfied by the assignment.
A clause is satisfied by an assignment $\alpha$ if it contains a literal that is satisfied by~$\alpha$.
Finally, a formula is satisfied by an assignment $\alpha$ if all its clauses are satisfied by $\alpha$.
A formula is \emph{satisfiable} if there exists an assignment that satisfies it and otherwise it is \emph{unsatisfiable}.
Two formulas are \emph{logically equivalent} if they are satisfied by the same assignments;
they are \emph{satisfiability equivalent} if they are either both satisfiable or both unsatisfiable.

\subsection{Clausal Proofs}

In the following, we introduce a formal notion of clause redundancy. 
%
A clause $C$ is \emph{redundant} with respect to a formula $F$ if $F$ and $F \land C$ are satisfiability equivalent.
%
For instance, the clause $C = x \lor y$ is redundant with respect to the formula $F = (\overline x \lor \overline y)$ since $F$
and $F \land C$ are satisfiability equivalent (although they are not logically equivalent).
This redundancy notion allows us to add redundant clauses to a formula without affecting its
satisfiability. 

Given a formula $F = \{C_1, \dots, C_m\}$, a \emph{clausal derivation} of a clause $C_n$ from $F$ is a sequence $C_{m+1}, \dots, C_n$ of clauses.
Such a sequence gives rise to \mbox{formulas} $F_m, F_{m+1}, \dots, F_n$, where $F_i = \{C_1, \dots, C_i\}$. We call $F_i$ the \emph{accumulated formula} corresponding to the \mbox{$i$-th} proof step.
A clausal derivation is \emph{correct} if every clause $C_i$ ($i > m$) is redundant with respect to the formula $F_{i-1}$ and if this redundancy can be checked in polynomial time with respect to the size of the proof.
A clausal derivation is a \emph{proof} of a formula $F$ if it derives the unsatisfiable empty clause. 
Clearly, since every clause-addition step preserves satisfiability, and since the empty clause
is always false, a proof of $F$ certifies the unsatisfiability of~$F$.
The proofs computed in this paper show that the chromatic a given graph is at least 5. 
We will also refer to proofs as {\em refutations} as they refute the existence of a valid 4-coloring.

\section{Clausal Proof Minimization}

SAT solving techniques are not only useful to validate the chromatic number of a graph, but they can also help reduce the size of the 
graph while preserving the chromatic number. The method works as follows. Given a graph $G$ with chromatic number $k$, first
generate the propositional formula $F$ that encodes whether the graph can be colored with $k-1$ colors. This formula is unsatisfiable.
Most SAT solvers can emit a proof of unsatisfiability. There exist several checkers for such proofs, even checkers that are
formally verified in the theorem provers {\sf ACL2}, {\sf Coq}, and {\sf Isabelle}~\cite{Cruz-Filipe2017,Lammich2017}. We used the
(unverified) 
checker {\sf DRAT-trim}~\cite{Heule:2013:trim} that allows minimizing the clausal proof as well as extracting an unsatisfiable core, i.e., a subformula that is also unsatisfiable.
From the unsatisfiable core one can easily extract a subgraph $G'$ of $G$ such that $G'$ also has chromatic 
number $k$. 

\subsection{Encoding}

We can compute the chromatic number of a graph $G$ as follows. Construct two formulas,  
one asking whether $G$ can be colored with $k-1$ colors, and one whether $G$ can be colored with $k$ colors. Now, 
$G$ has chromatic number $k$ if and only if the former is unsatisfiable while the latter is satisfiable. 

The construction of these two formulas can be achieved using the following encoding. Given a graph $G = (V,E)$
and a parameter $k$, the encoding uses $k|V|$ boolean variables $x_{v,c}$ with $v \in V$ and $c \in \{1,\dots,k\}$.
These variables have the following meaning: $x_{v,c}$ is true if and only if vertex $v$ has color $c$. Now we can encode
whether $G$ can be colored with $k$ colors:

\[
G_k : = \bigwedge_{v \in V} (x_{v,1} \lor \dots \lor x_{v,k}) \land \bigwedge_{\{v,w\} \in E} \bigwedge_{c \in \{1,\dots,k\}} (\overline x_{v,c} \lor \overline x_{w,c})
\]

The first type of clauses ensures that each vertex has at least one color, while the second type of clauses forces that two  
connected vertices are colored differently. Additionally, we could include clauses to require that each vertex has at most one color. However, these
clauses are redundant and would be eliminated by blocked clause elimination~\cite{BCE}, a SAT preprocessing technique. 

We added symmetry-breaking predicates~\cite{Crawford} during all experiments to speedup solving and proof minimization. 
The color-symmetries were broken by fixing the vertex at ($0$, $0$) to the first color, the vertex at ($1$, $0$) to the second color, and the vertex
at ($1/2$, $\sqrt{3}/2$) to the third color. These three points are at distance 1 from each other and occurred in all our graphs. The speedup is roughly a factor of
24 ($4 \cdot 3 \cdot 2)$, when trying to find a 4-coloring. 

We did not explore yet whether an encoding based on Zykov contraction~\cite{Zykov} would allow shorter proofs of unsatisfiability.
In essence, such an encoding would add variables and clauses that encode for a pair of vertices whether they have the same color. 
Solving graph coloring problems using such an extended encoding has been successful in the past~\cite{Schaafsma}.


\subsection{Graph Trimming}

Modern SAT solvers can emit clausal proofs. We used
the SAT solver {\sf Glucose}~\cite{glucose} to produce the proofs.
The most commonly supported 
format for clausal proofs is $\drat$, which computes the redundancy of clauses using the resolution asymmetric
tautology check~\cite{rules}. 
Some $\drat$ proof checkers can extract from a refutation an unsatisfiable core, i.e., a subformula
that is still unsatisfiable. When the formula expresses a graph coloring property, the unsatisfiable core represents
a subgraph with the same coloring property. The absence of the clause $(x_{v,1} \lor \dots \lor x_{v,k})$ in the core
shows that vertex $v$ can be removed, while the absence of all clauses $(\overline x_{v,c} \lor \overline x_{w,c})$
with $c \in \{1,\dots,k\}$ show that edge $\{v,w\}$ can be removed. When trying to find a small unit-distance graph with a given
chromatic number, we are interested in reducing the number of vertices. Although the proof checker can be easily modified to
ensure that no edges are removed, we achieved larger reductions by allowing edges to be deleted and then restoring edges
between vertices that survived the shrinking.

\subsection{Randomization}

SAT solvers and clausal-proof-minimization tools are deterministic. To increase the probability of finding small 
unit-distance graphs with chromatic number 5, we want to randomize the process and minimize many clausal
proofs. 

The proofs produced by SAT solvers depend heavily on the ordering of the clauses in the input file.
The initial heuristic ordering of the
variables is based on their occurrence in the input file. The earlier a variable occurs in the input file, the higher its
place in the ordering. Although more sophisticated initialization methods have been proposed, this method is effective in
practice. The effectiveness is caused by the typical encoding of a problem into propositional logic where
one starts with the more important variables. However, for our application there are no clear important variables.

Based on these observations, we applied the following lightweight randomization. First, we shuffle the input formula
and apply graph trimming on the result. When the clausal-proof-minimization tool is no longer able to remove vertices
from the graph, we shuffle the clauses of the current formula and produce a new clausal proof. Then we continue
graph trimming using the new formula and proof. This process is repeated until randomization cannot further 
reduce the size of the graph.


\subsection {Critical Graphs}

A graph is vertex/edge {\em critical} with respect to a given property if removing any vertex/edge would break that property.
Here we are interested in vertex critical graphs with respect to the chromatic number. 
Graph trimming as described above would remove most redundant vertices of the graph and the randomization method 
allows shrinking the graph even further. However, in most cases the reduced graphs are not critical: There still exist some
vertices that can be removed while preserving the chromatic number.

Both the SAT solver and the clausal-proof-minimization tool aim to find a relatively short argument (i.e., clausal proof) explaining why
imply the fewest number of involved vertices. In fact, an argument can frequently be shortened 
by using redundant (non-critical) vertices. 

In the final step we therefore make the graph critical by the following procedure. Randomly pick a vertex from the 
graph and determine the chromatic number of the graph without it. If the chromatic number is not changed, then the
vertex is removed from the graph. This process is repeated until all remaining vertices have been determined to be
critical.

Instead of using this naive method to make the graph critical, we could have used more sophisticated tools
that compute a minimal unsatisfiable core from the propositional formula. However, these tools did not improve 
the performance or the size in an observable way.


\subsection{Validation}

Determining whether a set of points form a unit-distance graph with chromatic number 5 requires two
checks: (i) does the corresponding graph have chromatic number 5; and (ii) is the distance between two connected points
exactly 1. The techniques discussed in this paper can easily perform the first check. SAT solvers can compute valid 5-coloring
for the critical graphs in a fraction of a second. The proofs showing that there exists no 4-coloring are actually quite small:
between $14\,000$ and $19\,000$ clause addition steps. Proofs of that size
can be checked in roughly a second even with formally verified checkers. We used the {\sf DRAT-trim} tool~\cite{Heule:2013:trim} to validate them.
Proofs of recently solved hard-combinatorial problems, 
such as the Pythagorean Triples and Schur Number Five, are much larger: roughly 1 trillion and 10 trillion clause addition steps, respectively~\cite{ptn,S5}.

For the second check we used a tool based on Gr\"obner basis, available at \url{http://fmv.jku.at/dist1sqrtgb/}, to validate for every edge in the graph that 
the corresponding points are exactly 1 apart. The tool produces files that can be validated using 
{\sf Singular}~\cite{Singular} and {\sf pactrim}~\cite{RitircBiereKauers-SCSC18}. There is no need to check whether all edges are present as missing edges can only decrease the chromatic number. 
Checking only the correctness of the edges in the graph is cheap. The total validation time for our smallest
critical graphs is about a second or two.

\begin{figure}[b]
\includegraphics[width=0.45\textwidth]{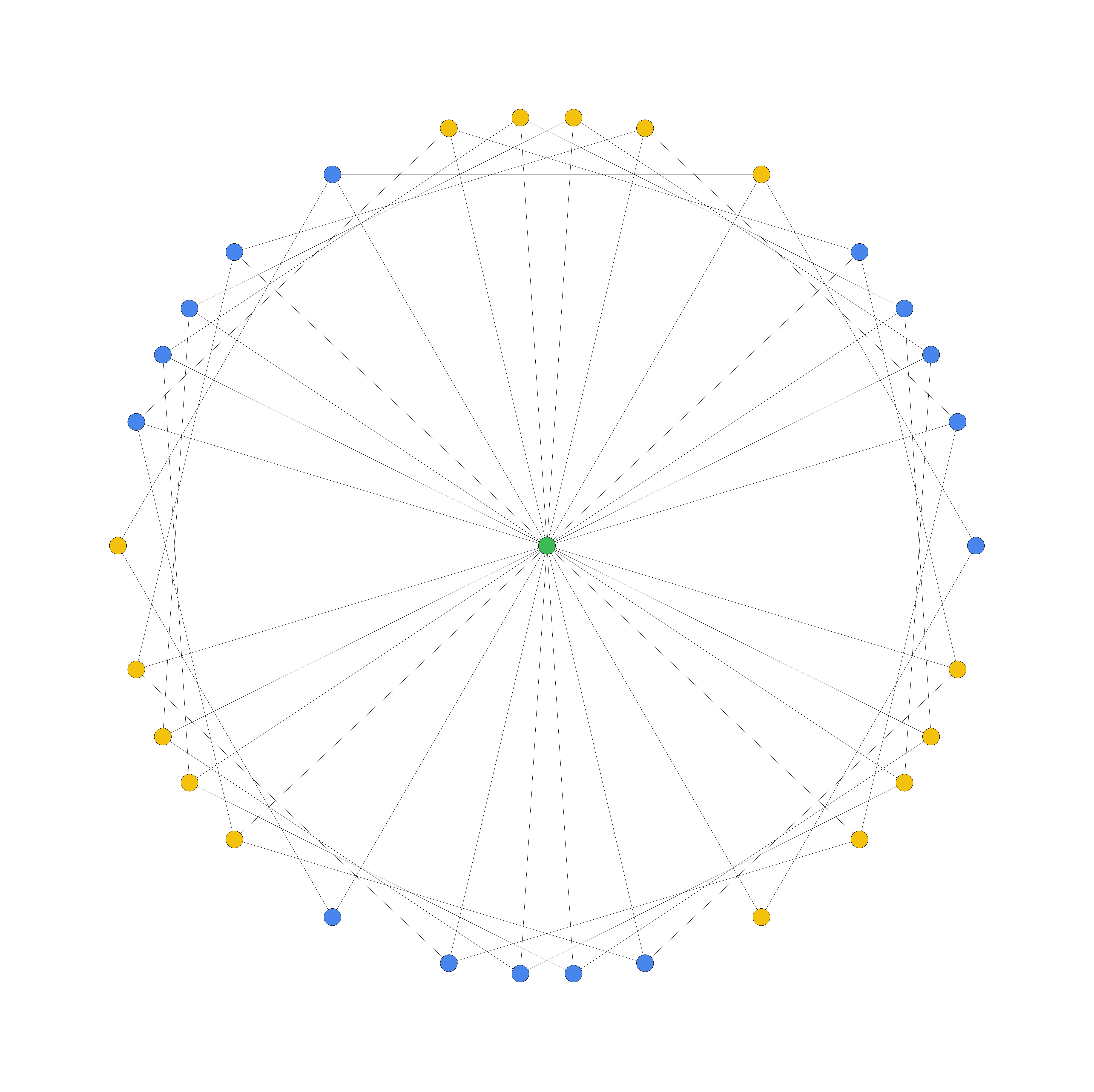}
\hfill
\includegraphics[width=0.45\textwidth]{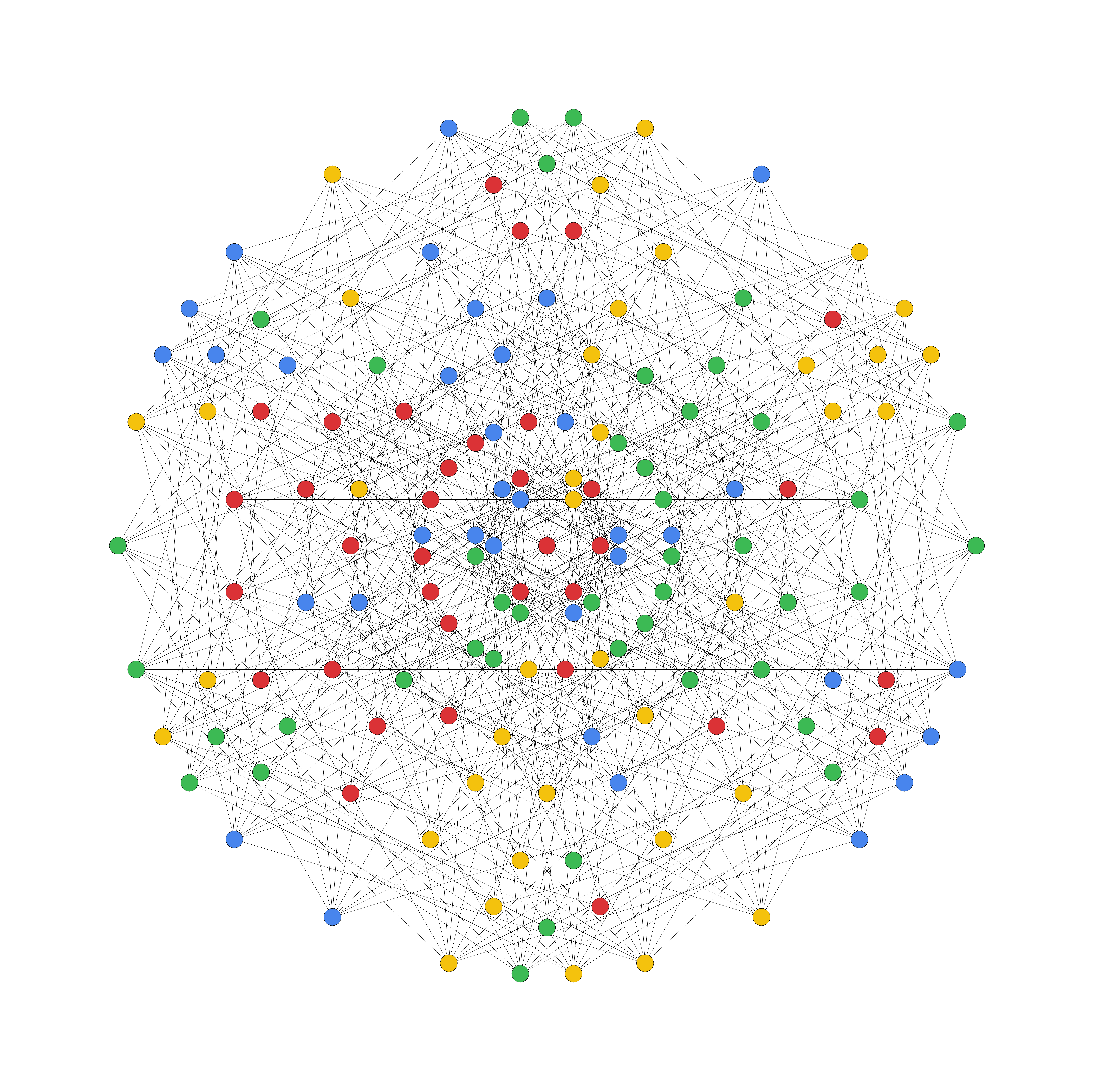}
\caption{Left, a 3-coloring of de Grey's graph $V_{31}$. Right, a 4-coloring of $V_{151}$ being $V_{31} \oplus V_{31}$ without the
vertices more than unit distance apart from the center.}
\label{fig:V}
\end{figure}

\section{Results}

In this section we discuss the various techniques that we used and developed to obtain small unit-distance graphs 
with chromatic number 5. The techniques were originally designed for verification purposes and applying them to graph minimization
is novel and unexpected. The main strategy is to start with a large graph and shrink it using clausal proof minimization.
We minimized several large graphs with various heuristics and most reduced graphs consisted of 800 to 900 vertices. However,
we were able to produce several graphs of 553 vertices using three techniques. The first technique (Section~\ref{sec:small})
enabled producing graphs with less than 700 vertices consistently. Second, we obtained graphs with just over 600 vertices
by shrinking merged copies of graphs with less than 700 vertices (Section~\ref{sec:merge}). Finally, we added some points far away
from the origin in order to eliminate more points close to the origin (Section~\ref{sec:outer}). The graphs and corresponding proofs mentioned in this section are
available at \url{https://github.com/marijnheule/CNP-SAT}.

\subsection{Finding a small symmetric subgraph}
\label{sec:small}

The main building block of our graphs is de Grey's $V_{31}$~\cite{DeGrey}: five 7-wheels with a common central vertex. The graph 
$V_{31}$ has 31 vertices, 60 edges, is 3-colorable, and all points are in the field $\mathbb{Q}[\sqrt{3}, \sqrt{11}]$. These points can 
be obtain by applying $\theta_1^{j}\theta_3^{k}$ on point ($1$, $0$) around ($0$, $0$) with $j \in \{0,1,2,3,4,5\}$ and $k \in \{-1,-\frac{1}{2},0,\frac{1}{2},1\}$.
A visualization of this graph is shown in Figure~\ref{fig:V} (left). During an early stage of the experimentation,
we observed that the graph $(V_{31} \oplus V_{31} \oplus V_{31}) \cup \theta_4(V_{31} \oplus V_{31} \oplus V_{31})$ has chromatic number 5. 
Furthermore, all points that are further away than 2 of the center can be removed without affecting the chromatic number. 

\begin{figure}[b]
\includegraphics[width=.95\textwidth]{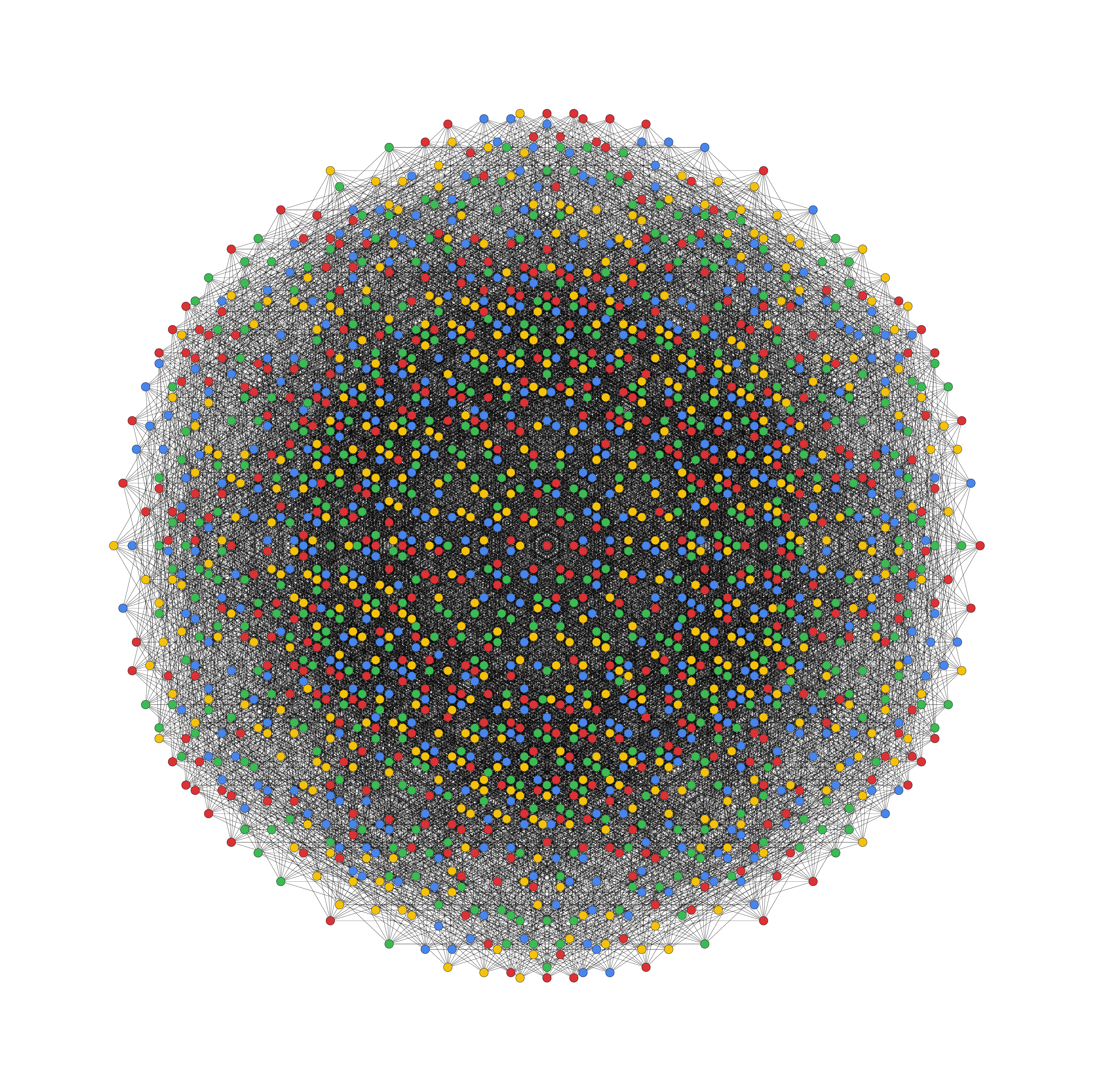}
\caption{A 4-coloring of $V_{1939}$, which is the Minkowski sum of $V_{31}$ and $V_{151}$.}
\label{fig:V1939}
\end{figure}

Instead of removing the points at distance larger than 2 from the center, we constructed the following graph. Let $V_{151}$ be the
Minkowski sum of $V_{31}$ and $V_{31}$ without the points at distance larger than 1 from the center. This graph has 151 vertices and
510 edges and is shown in Figure~\ref{fig:V} (right). Now let $V_{1939}$ be the Minkowski sum of $V_{31}$ and $V_{151}$. This graph 
is shown in Figure~\ref{fig:V1939}. The graph $V_{1939} \cup \theta_4(V_{1939})$ has chromatic number 5 as well.

We applied clausal proof minimization on the 
formula that encodes whether the graph $V_{1939} \cup \theta_4(V_{1939})$ is 4-colorable. 
Most random probes of clausal proof minimization produced a subgraph of $V_{1939} \cup \theta_4(V_{1939})$ with slightly more than 800 vertices.
Occasionally it produces graphs with fewer than 700 vertices, while never producing graphs in the range of 700 to 800 vertices. 

\begin{figure}[b]
\centering
\includegraphics[width=0.43\textwidth]{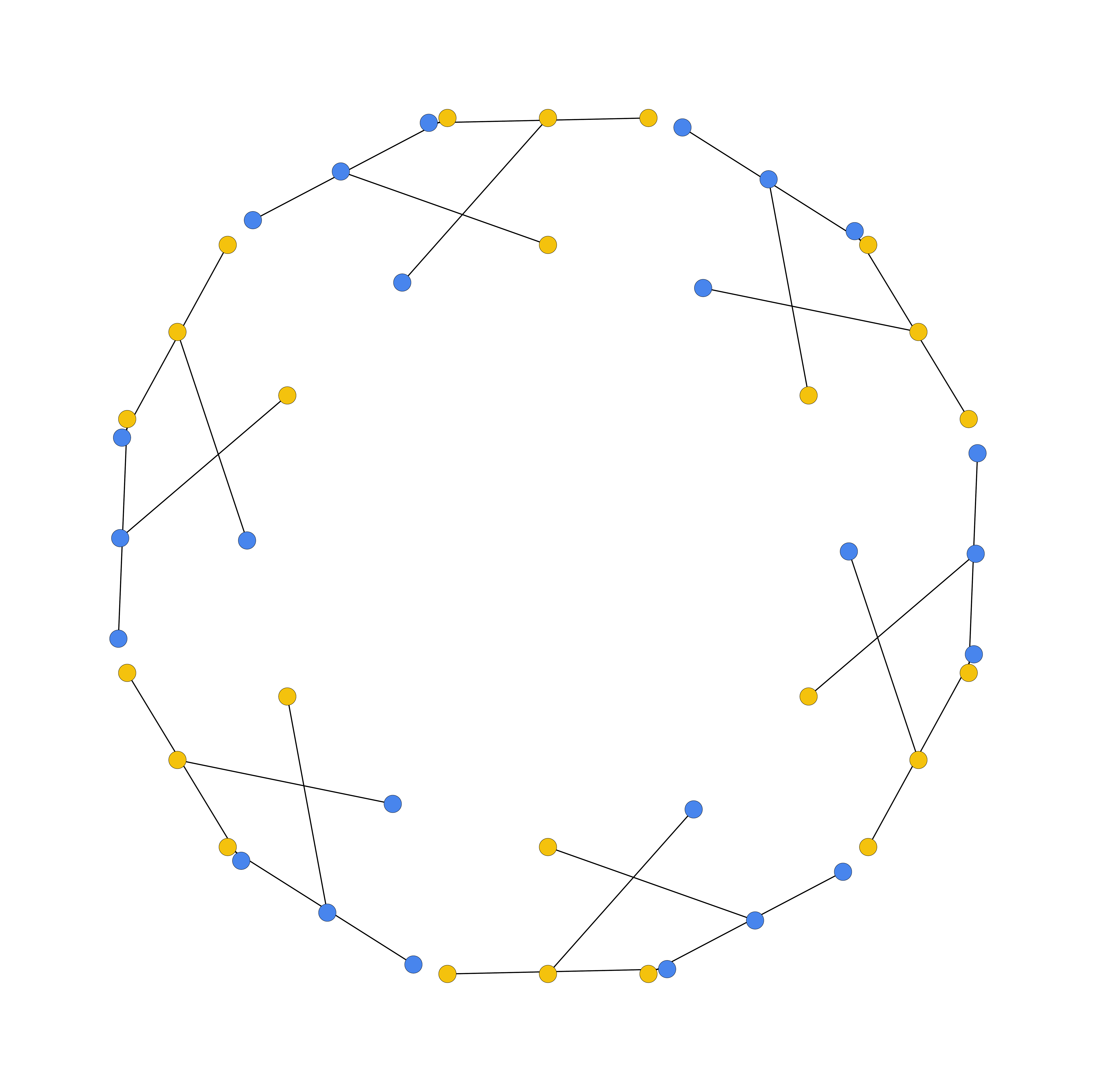}
\caption{A visualization of the edges between $S_{199} \cup \theta_4(S_{199})$ and the involved vertices. Vertices originating from
$S_{199}$ and $\theta_4(S_{199})$ are colored yellow and blue, respectively.}
\label{fig:T4S199}
\end{figure}

Closer examination of the minimized graphs with fewer than 700 vertices revealed that only a small fraction of the points (always less than 200 vertices) are in the
field $\mathbb{Q}[\sqrt{3}, \sqrt{11}]$. These points originate from the subgraph $V_{1939}$, while the other points originate from the subgraph $\theta_4(V_{1939})$.
Other patterns can be observed in the graphs with fewer than 700 vertices: there were at least 12 points in the field $\mathbb{Q}[\sqrt{3}, \sqrt{11}]$ at distance 2,
while the graphs with more than 800 vertices had fewer than three such points. Hence keeping the points at distance 2 appears crucial to find smaller graphs.
Rotation $\theta_4$ does not only add edges between points at distance 2 (by construction), but also between points at other distances. 
In fact, half the edges between points in $V_{1939}$
and $\theta_4(V_{1939})$ are due to points that are closer to the center: i.e.,  at $\frac{\sqrt{33} + 1}{2\sqrt{3}}$ and $\frac{\sqrt{33} - 1}{2\sqrt{3}}$ from the origin.
Figure~\ref{fig:T4S199} shows the newly introduced edges due to $\theta_4$.


\begin{figure}[t]
\includegraphics[width=0.95\textwidth]{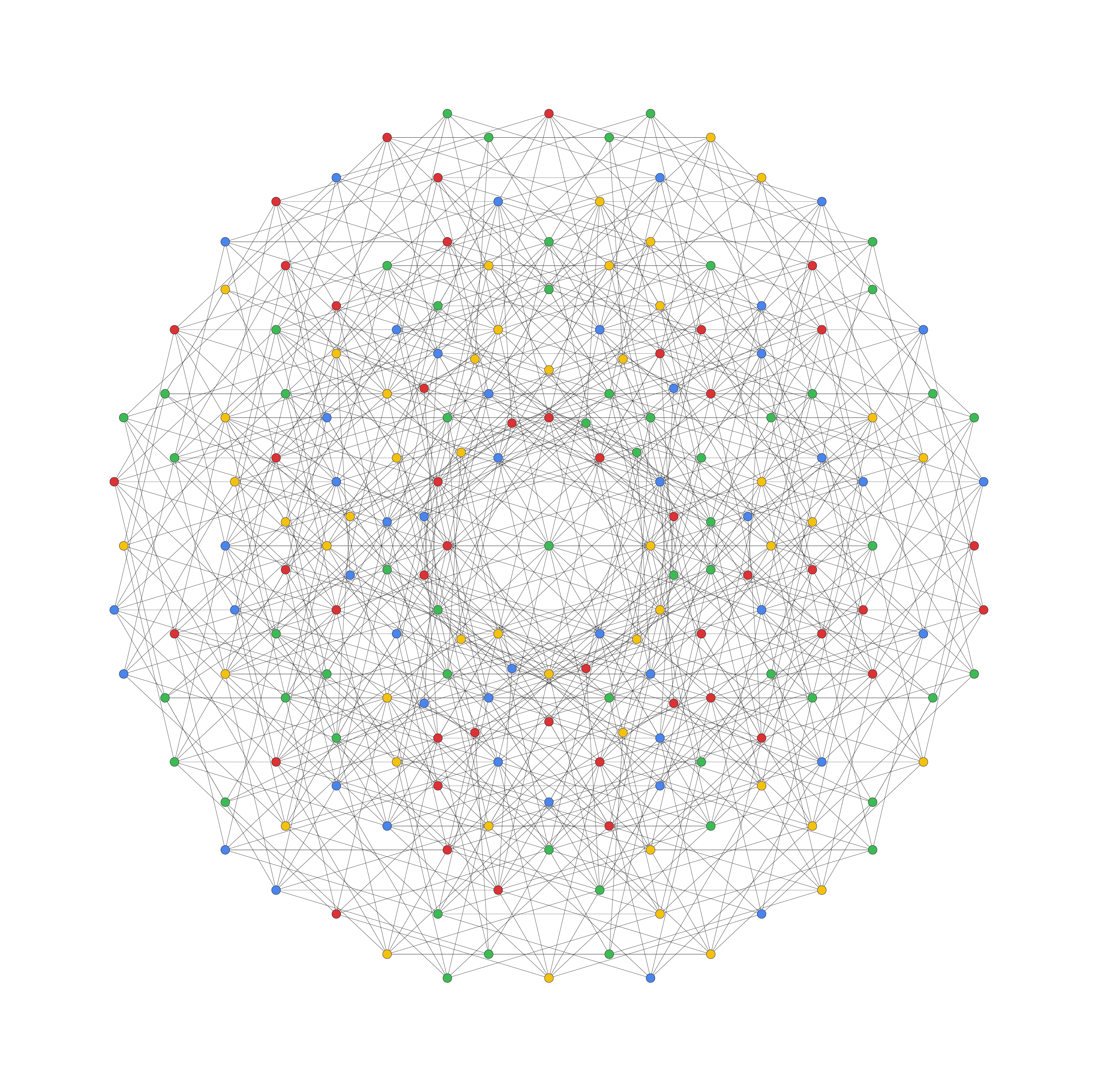}
\caption{A 4-coloring of the graph $S_{199}$, a symmetric subgraph of $V_{1939}$, that occurred in several early records.}
\label{fig:S199}
\end{figure}

Visualizing the points in the field $\mathbb{Q}[\sqrt{3}, \sqrt{11}]$ reveals that they are highly symmetric: both reflection in the horizontal axis and a rotation of $\theta_1 = 60^{\circ}$
map the points onto themselves. Figure~\ref{fig:S199} shows this visualization.
Shown is a 199-vertex graph with 888 edges at unit distance, which we call $S_{199}$. The minimized graphs did not fully produce $S_{199}$, but always yielded
a subgraph that missed a handful (up to a dozen) vertices in various locations. There exist many 4-colorings of $S_{199}$, but we observed no clear pattern.

\begin{figure}[t]
\includegraphics[width=\textwidth]{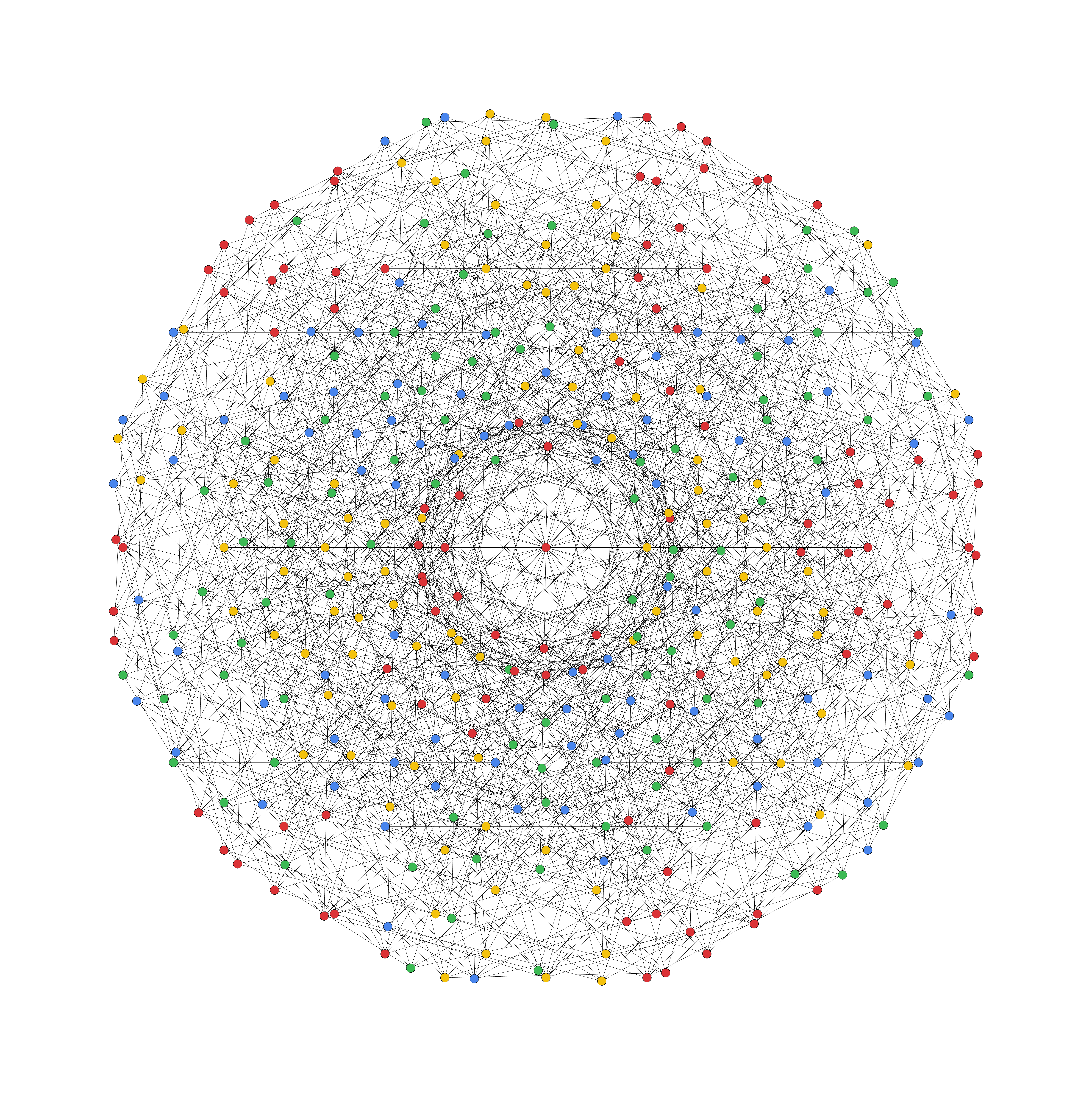}
\caption{A 4-coloring of the graph $S_{199} \cup \theta_4(S_{199})$.}
\label{fig:T7plot}
\end{figure}

Interesting patterns emerge when merging $S_{199}$ and $\theta_4(S_{199})$, as shown in Figure~\ref{fig:T7plot}. Notice that points that are close to each other
frequently have the same color. More importantly, roughly half of the vertices that are close to distance 2 from the center have the same color as the central vertex. 
In later experiments we minimized the graph $V_{1939} \cup \theta_4(S_{199})$, which allowed us to consistently produce unit-distance graphs with fewer than 700 vertices.
We suspect that the above mentioned patterns contribute to the lack of a 4-coloring of $V_{1939} \cup \theta_4(S_{199})$. Notice that $V_{1939}$ has $S_{199}$ as a subgraph.

\subsection {Merging Critical Graphs}
\label{sec:merge}

In order to further produce smaller unit-distance graphs with chromatic number 5, we selected two critical graphs obtained earlier, merged 
them, and applied clausal proof minimization again. There are a significant number of options to merge two graphs and we experimented with a variety 
of these. The most effective merging strategy in our experiments turned out to be rotating the graphs along the central vertex in
such a way that a vertex in one graph at unit distance from the center is merged with a vertex from the other graph at unit distance.
Although two different critical graphs can be used for merging, we observed that it is also effective to merge two copies of the same
critical graph. 

The minimization procedure frequently produced a graph that was larger than both the critical graphs that were merged. Only three kinds of 
rotations occasionally resulted in smaller graphs. The first two rotations are $\theta_1^k$ and $\theta_3^k$ for a small value of $k$. 
These rotations clearly increase the average vertex degree by merging vertices and introducing edges between points at distance 1 ($\theta_1^k$) and
$\sqrt{3}$ ($\theta_3^k$) from the origin. Most other rotations result in little interaction between the two critical graphs and thus
hardly increase the average vertex degree.

\begin{figure}
\includegraphics[width=\textwidth]{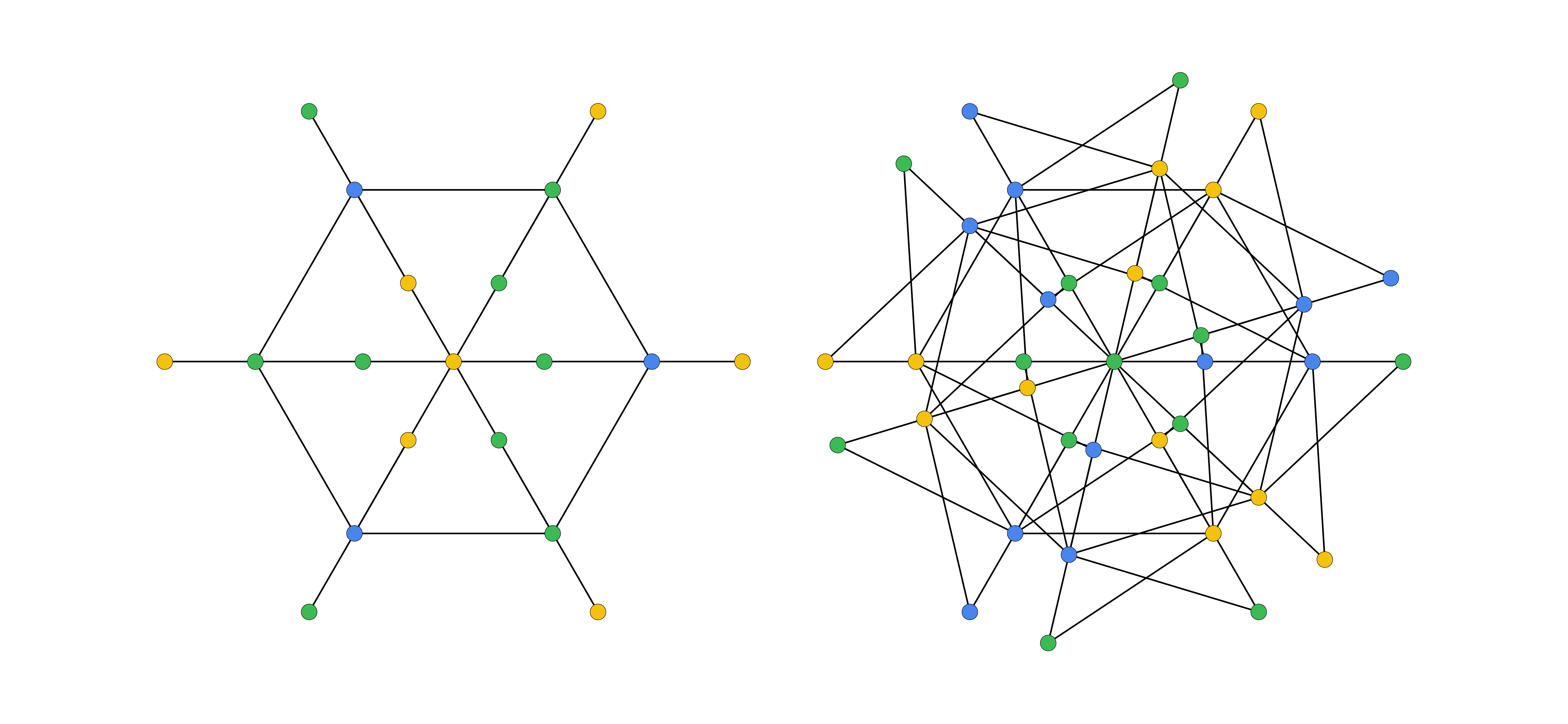}
\caption{A rotation by $\theta_3^{\frac{1}{2}}$ connects points at different distances from the origin.
Left, a subgraph of $V_{31} \oplus V_{31}$ consisting of three 7-wheels with radii $1$,  $\frac{\sqrt{33} + 3}{6}$, and $\frac{\sqrt{33} - 3}{6}$. 
Right, two copies of that graph with one copy turned by this rotation.}
\label{fig:3wheel}
\end{figure}

The most effective rotations introduce edges between the points at distance $1$ and the distances 
$\frac{\sqrt{33} + 3}{6}$ and $\frac{\sqrt{33} - 3}{6}$ from the origin, thereby increasing the average vertex degree significantly. Figure~\ref{fig:3wheel} illustrates this by showing three
7-wheel graphs with the radii $1$,  $\frac{\sqrt{33} + 3}{6}$, and $\frac{\sqrt{33} - 3}{6}$ (left) and two copies of this graph
rotated in such a way that the points on these distances become connected (right). The graph on the left has average vertex degree
$\frac {36}{19}$, while the graph of the right has average vertex degree $\frac{120}{37}$.
A rotation by $\theta_3^{\frac{1}{2}}$ for example achieves this and maps
point $(1,0)$ onto point $(\frac{\sqrt{33}}{6}, \frac{\sqrt{3}}{6})$.
Both points are at unit distance from the origin
and both are part of $V_{31}$ and of most other graphs that we used in the experiments. 

\begin{figure}[b]
\includegraphics[width=0.95\textwidth]{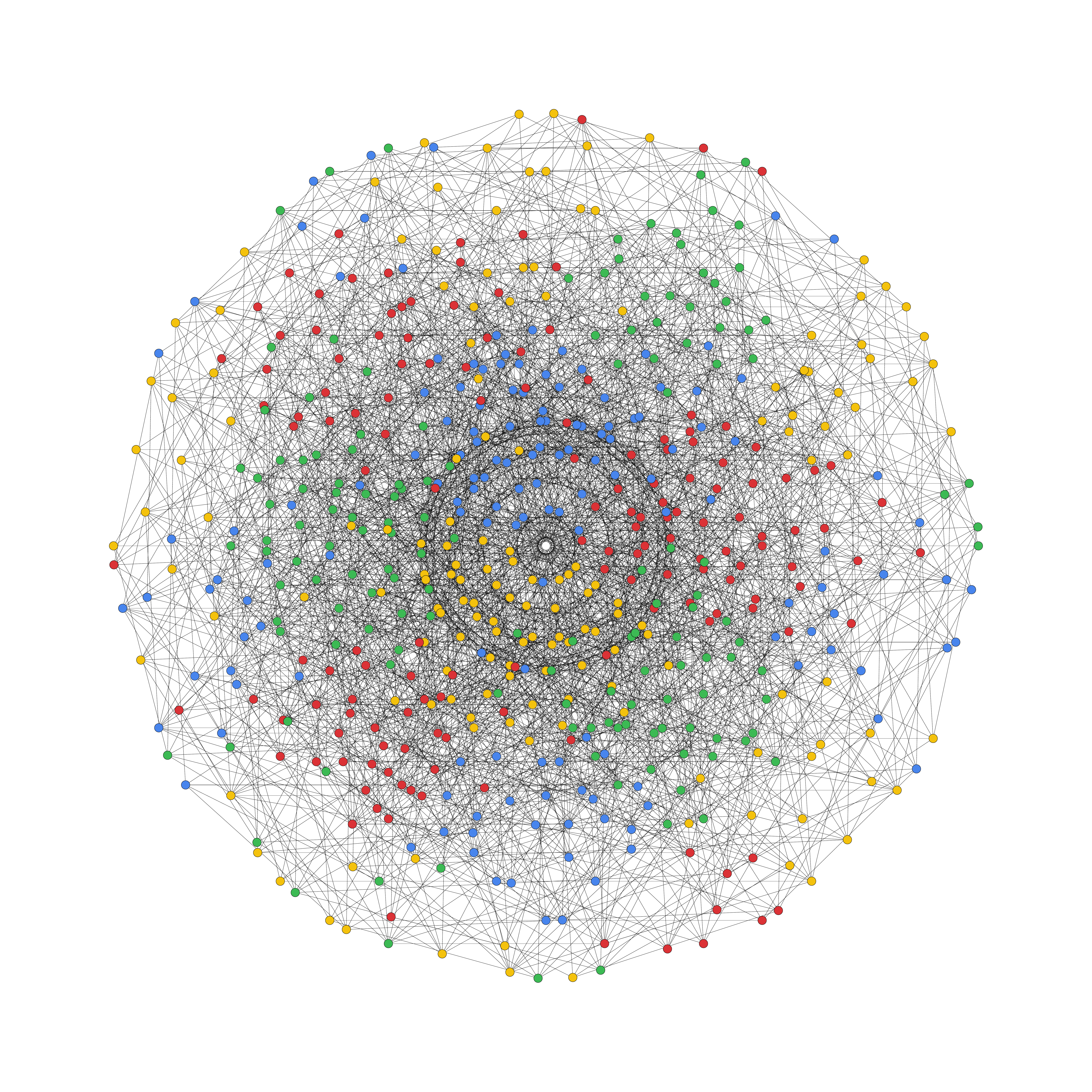}
\caption{A visualization of a 610-vertex unit-distance graph with chromatic number 5. 
Five colors are used for the vertices. Only the center uses the fifth color (white).}
\label{fig:610}
\end{figure}

The combination of merging and minimization only introduced vertices in the field $\mathbb{Q}[\sqrt{3}, \sqrt{11}]$. The smallest graphs contained 
roughly 50 vertices that do not occur in $(V_{31} \oplus V_{31} \oplus V_{31})$ and thus not in $V_{1939}$.

The smallest graph that we found using the techniques discussed so far contains 610 vertices and 3000 edges.
This graph is shown in Figure~\ref{fig:610} using a 5-coloring in which only the central vertex has the fifth color.
Recall that this graph is vertex critical. 
Hence our graph possesses such a coloring in which any vertex
can be the only one with the fifth color.

\subsection {Minimizing the Small Part}
\label{sec:outer}

The critical graphs found so far can be partitioned into two parts:  a subgraph of $\theta_4(S_{199})$ and the 
subgraph induced by the remaining vertices.
We refer to the former as the small part, as it consists typically of only 187 vertices, and to the latter as the large part.
In all statements regarding the size of these graphs, we count the central vertex in both parts.
All points in the smallest critical graphs are at distance 2 or less from the center. 
Several approaches have been examined in order to find unit-distance graphs with fewer than 600 vertices. Only 
one approach was effective.

\begin{figure}[h]
\includegraphics[width=0.95\textwidth]{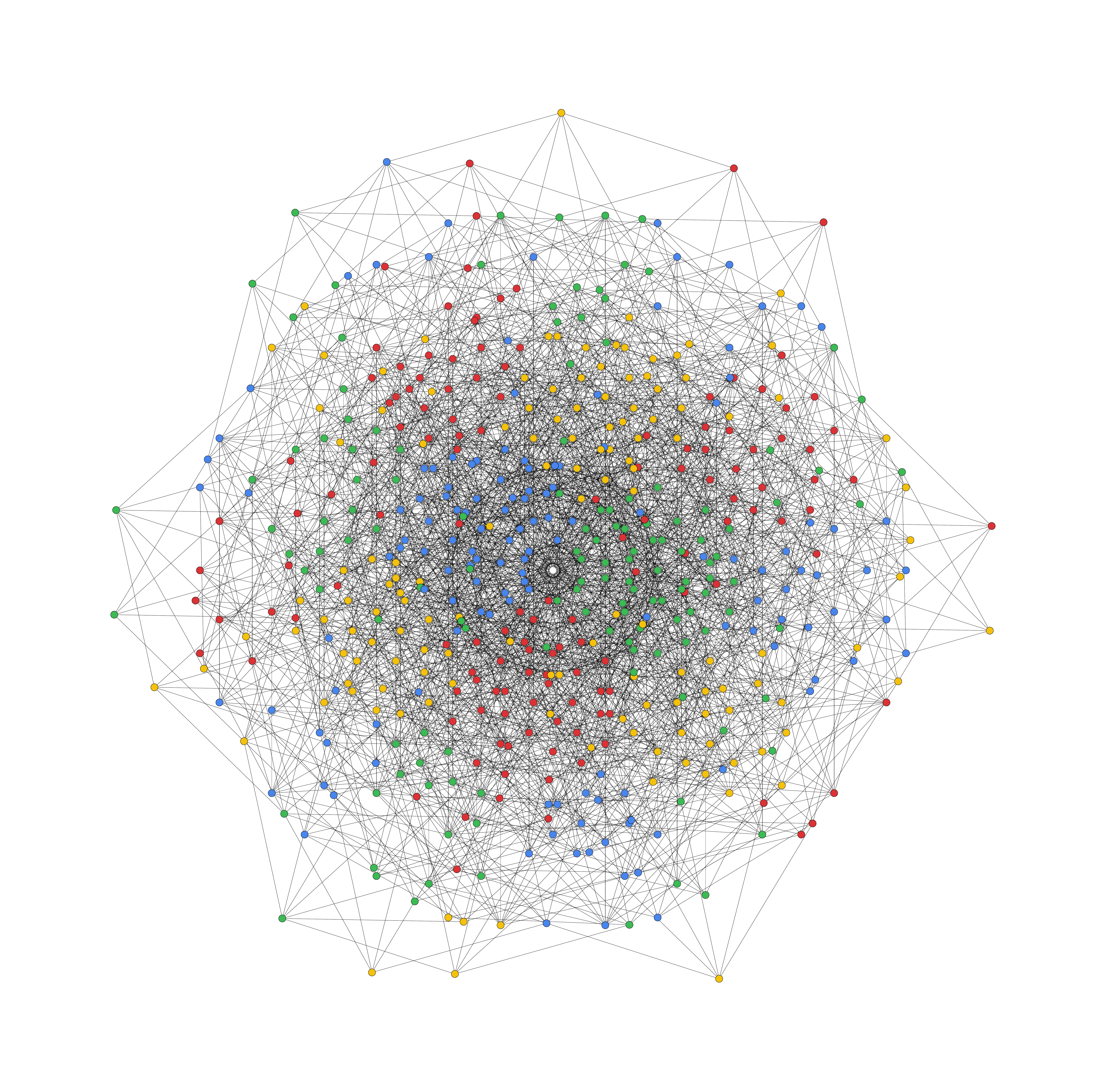}
\caption{A visualization of a 553-vertex unit-distance graph with chromatic number 5. 
Five colors are used for the vertices. Only the center uses the fifth color (white).}
\label{fig:553}
\end{figure}

We focussed on adding points that are further away than 2 from the origin in order to remove more inner vertices.
Adding points from the field $\mathbb{Q}[\sqrt{3}, \sqrt{11}]$ may allow reducing the large part, but none of the experiments were successful. 
However, we were able to substantially reduce the small part using this strategy. The most effective approach was as follows.
We first constructed the Minkowski sum of $\theta_4(S_{199})$ and $\theta_4(S_{199})$ and removed all points that
were less or equal than 2 away from the origin. This graph consists of 2028 vertices. All points were added to the smallest critical graphs
that were found in the earlier steps, followed by clausal proof minimization. This resulted in a dozen graphs with 553 vertices and
(on average) 2720 edges. Figure~\ref{fig:553} shows one of these graphs, which we refer to as $G_{553}$. 
Practically all vertices that were removed during minimization originated from the small part. This 
part was reduced to 133 or 134 vertices.  The 553-vertex graphs appear less symmetric compared to the earlier graphs. This is caused by the few vertices that are
far from the origin.

The unit-distance graphs with 553 vertices are vertex critical, but not edge critical. The proofs of unsatisfiability show that many of 
the edge clauses can be removed without introducing a 4-coloring. Randomly removing edges until fixpoint eliminates about 270 edges
(close to $10\%$) of these graphs.

Remarkably, all critical graphs have a handful of vertices with degree 4. If we removed such a vertex from the 
graph, all its four neighbors would have a different color in all valid 4-colorings. Graph $S_{199}$ has even 12 vertices with degree 4.
Reducing a 553-vertex graph to become edge critical will increase
the number vertices with degree 4 to roughly 12. These vertices tend to be evenly distributed between the small and large parts of the graph.

\subsection{Analysis}

The sizes of the small and the large parts of the 553-vertex graphs suggest that they play a different role in eliminating 4-colorings. We analyzed
the 4-colorings of both parts when restricted to the key vertices, i.e, the ones that connect the parts. These are the vertices at distance
$\frac{\sqrt{33} - 1}{2\sqrt{3}}$, $\frac{\sqrt{33} + 1}{2\sqrt{3}}$, and 2 from the origin. Recall that Figure~\ref{fig:T4S199} shows the 
interaction between these vertices. The 553-vertex graphs have all 24 vertices occurring in $S_{199} \cup \theta_4(S_{199})$ at distance 2 from the
original and most of the vertices at distance $\frac{\sqrt{33} - 1}{2\sqrt{3}}$, $\frac{\sqrt{33} + 1}{2\sqrt{3}}$ from the origin.

\begin{figure}[t]
\centering
\includegraphics[width=0.95\textwidth]{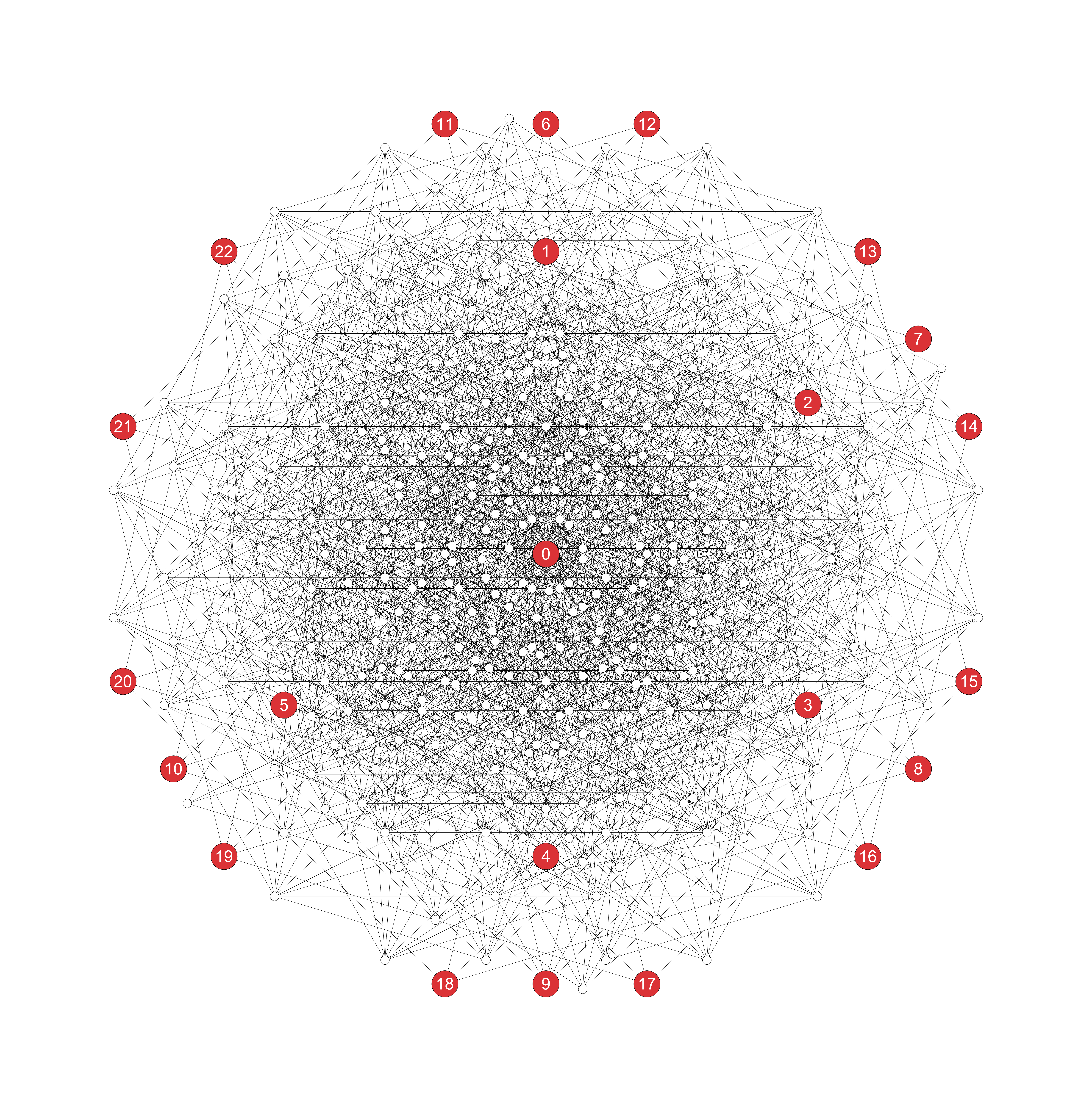}
\caption{The 420-vertex large part of $G_{553}$ with the key vertices marked with numbers.}
\label{fig:numberL}
\end{figure}

Figure~\ref{fig:numberL} shows the large part of $G_{553}$ in which the central and key vertices  are numbered. The other vertices in the large
part significantly restrict the number of different 4-colorings of these key vertices. In fact, there are only twenty 4-colorings of these vertices such
that they either have the same color as the central vertex or a different color. Table~\ref{tab:solutions} shows the details. The clearest pattern
is that either the vertices $v_1$ to $v_5$ (at distance $\frac{\sqrt{33} - 1}{2\sqrt{3}}$) have the same color as the central vertex
or the vertices $v_6$ to $v_{10}$ (at distance $\frac{\sqrt{33} + 1}{2\sqrt{3}}$) have the same color as the central vertex. This pattern heavily 
constrains the vertices in the small part at those distances: Either all vertices at distance $\frac{\sqrt{33} + 1}{2\sqrt{3}}$ or all vertices at distance
$\frac{\sqrt{33} - 1}{2\sqrt{3}}$ from the origin must have a different color than the central vertex. Other patterns can be observed 
as well. For example, if the missing two vertices would be added, then most colorings can be obtained by a rotation of 60 degrees. 

In contrast, the small part hardly constrains its key vertices since as many as $9974$ 4-colorings are allowed. However, none of these many 4-colorings
is compatible with the twenty 4-colorings of the large part. Enforcing the above mentioned pattern, i.e., either all vertices at distance
$\frac{\sqrt{33} + 1}{2\sqrt{3}}$ or all vertices at distance $\frac{\sqrt{33} - 1}{2\sqrt{3}}$ from the origin must have a different color than the central vertex,
reduces the number of 4-colorings of the small part to 1353. 


\newcolumntype{P}[1]{>{\centering\arraybackslash}p{#1}}

\begin{table}[t]
\caption{A list of the twenty 4-colorings of $v_1$ to $v_{22}$ in $G_{553}$. The $1$s denote that the vertex has the same color as the central vertex, while
the $0$s denote that the vertex is colored differently than the central vertex. The bold $1$s indicate the main pattern.}
\label{tab:solutions}
\centering
\taburulecolor{lightgray}
\begin{tabular}{|@{}P{0.55cm}@{}|@{}P{0.55cm}@{}|@{}P{0.55cm}@{}|@{}P{0.55cm}@{}|@{}P{0.55cm}@{}||@{}
                          P{0.55cm}@{}|@{}P{0.55cm}@{}|@{}P{0.55cm}@{}|@{}P{0.55cm}@{}|@{}P{0.55cm}@{}||@{}P{0.55cm}@{}|@{}
                          P{0.55cm}@{}|@{}P{0.55cm}@{}|@{}P{0.55cm}@{}|@{}P{0.55cm}@{}|@{}P{0.55cm}@{}|@{}P{0.55cm}@{}|@{}P{0.55cm}@{}|@{}P{0.55cm}@{}|@{}P{0.55cm}@{}|@{}P{0.55cm}@{}|@{}P{0.55cm}@{}|}
\hhline{-----||-----||------------}
$v_1$ & $v_2$ &$v_3$ &$v_4$ &$v_5$ &$v_6$ &$v_7$ &$v_8$ &$v_9$ &$v_{10}$ &$v_{11}$ &
$v_{12}$ &$v_{13}$ &$v_{14}$ &$v_{15}$ &$v_{16}$ &$v_{17}$ &$v_{18}$ &$v_{19}$ &$v_{20}$ &$v_{21}$ &$v_{22}$ \\ \hhline{-----||-----||------------}\noalign{\vspace{3pt}}\hhline{-----||-----||------------}
0 & 0 & 0 & 0 & 0 &  {\bf 1} & {\bf 1} & {\bf 1} & {\bf 1} & {\bf 1} & 1 & 1 & 1 & 1 & 1 & 1 & 1 & 1 & 1 & 1 & 1 & 1 \\ \hhline{-----||-----||------------}
0 & 0 & 0 & 0 & 0 &  {\bf 1} & {\bf 1} & {\bf 1} & {\bf 1} & {\bf 1} & 0 & 1 & 0 & 1 & 1 & 1 & 1 & 0 & 1 & 0 & 1 & 1 \\ \hhline{-----||-----||------------}
1 & 0 & 0 & 0 & 0 &  {\bf 1} & {\bf 1} & {\bf 1} & {\bf 1} & {\bf 1} & 1 & 1 & 0 & 1 & 0 & 1 & 1 & 1 & 1 & 0 & 1 & 0 \\ \hhline{-----||-----||------------}
0 & 1 & 0 & 0 & 0 &  {\bf 1} & {\bf 1} & {\bf 1} & {\bf 1} & {\bf 1} & 1 & 0 & 1 & 1 & 0 & 1 & 0 & 1 & 1 & 1 & 1 & 0 \\ \hhline{-----||-----||------------}
0 & 0 & 1 & 0 & 0 &  {\bf 1} & {\bf 1} & {\bf 1} & {\bf 1} & {\bf 1} & 1 & 0 & 1 & 0 & 1 & 1 & 0 & 1 & 0 & 1 & 1 & 1 \\ \hhline{-----||-----||------------}
0 & 0 & 0 & 1 & 0 &  {\bf 1} & {\bf 1} & {\bf 1} & {\bf 1} & {\bf 1} & 1 & 1 & 1 & 0 & 1 & 0 & 1 & 1 & 0 & 1 & 0 & 1 \\ \hhline{-----||-----||------------}
0 & 0 & 0 & 0 & 1 &  {\bf 1} & {\bf 1} & {\bf 1} & {\bf 1} & {\bf 1} & 0 & 1 & 1 & 1 & 1 & 0 & 1 & 0 & 1 & 1 & 0 & 1 \\ \hhline{-----||-----||------------}
1 & 0 & 0 & 1 & 0 &  {\bf 1} & {\bf 1} & {\bf 1} & {\bf 1} & {\bf 1} & 1 & 1 & 0 & 0 & 0 & 0 & 1 & 1 & 0 & 0 & 0 & 0 \\ \hhline{-----||-----||------------}
0 & 1 & 0 & 0 & 1 &  {\bf 1} & {\bf 1} & {\bf 1} & {\bf 1} & {\bf 1} & 0 & 0 & 1 & 1 & 0 & 0 & 0 & 0 & 1 & 1 & 0 & 0 \\ \hhline{-----||-----||------------}
0 & 0 & 1 & 0 & 0 &  {\bf 1} & {\bf 1} & {\bf 1} & {\bf 1} & {\bf 1} & 0 & 0 & 0 & 0 & 1 & 1 & 0 & 0 & 0 & 0 & 1 & 1 \\ \hhline{-----||-----||------------} \noalign{\vspace{3pt}}\hhline{-----||-----||------------}
{\bf 1} & {\bf 1} & {\bf 1} & {\bf 1} & {\bf 1} & 0 & 0 & 0 & 0 & 0 & 1 & 1 & 1 & 1 & 1 & 1 & 1 & 1 & 1 & 1 & 1 & 1 \\ \hhline{-----||-----||------------}
{\bf 1} & {\bf 1} & {\bf 1} & {\bf 1} & {\bf 1} & 0 & 0 & 0 & 0 & 0 & 1 & 0 & 1 & 1 & 1 & 1 & 1 & 1 & 0 & 1 & 0 & 0 \\ \hhline{-----||-----||------------}
{\bf 1} & {\bf 1} & {\bf 1} & {\bf 1} & {\bf 1} & 1 & 0 & 0 & 0 & 0 & 0 & 0 & 1 & 0 & 1 & 1 & 1 & 1 & 1 & 1 & 0 & 1 \\ \hhline{-----||-----||------------}
{\bf 1} & {\bf 1} & {\bf 1} & {\bf 1} & {\bf 1} & 0 & 1 & 0 & 0 & 0 & 0 & 1 & 0 & 0 & 1 & 0 & 1 & 1 & 1 & 1 & 1 & 1 \\ \hhline{-----||-----||------------}
{\bf 1} & {\bf 1} & {\bf 1} & {\bf 1} & {\bf 1} & 0 & 0 & 1 & 0 & 0 & 1 & 1 & 0 & 1 & 0 & 0 & 1 & 0 & 1 & 1 & 1 & 1 \\ \hhline{-----||-----||------------}
{\bf 1} & {\bf 1} & {\bf 1} & {\bf 1} & {\bf 1} & 0 & 0 & 0 & 1 & 0 & 1 & 1 & 1 & 1 & 0 & 1 & 0 & 0 & 1 & 0 & 1 & 1 \\ \hhline{-----||-----||------------}
{\bf 1} & {\bf 1} & {\bf 1} & {\bf 1} & {\bf 1} & 0 & 0 & 0 & 0 & 1 & 1 & 1 & 1 & 1 & 1 & 1 & 0 & 1 & 0 & 0 & 1 & 0 \\ \hhline{-----||-----||------------}
{\bf 1} & {\bf 1} & {\bf 1} & {\bf 1} & {\bf 1} & 1 & 0 & 0 & 1 & 0 & 0 & 0 & 1 & 0 & 0 & 1 & 0 & 0 & 1 & 0 & 0 & 1 \\ \hhline{-----||-----||------------}
{\bf 1} & {\bf 1} & {\bf 1} & {\bf 1} & {\bf 1} & 0 & 1 & 0 & 0 & 1 & 0 & 1 & 0 & 0 & 1 & 0 & 0 & 1 & 0 & 0 & 1 & 0 \\ \hhline{-----||-----||------------} 
{\bf 1} & {\bf 1} & {\bf 1} & {\bf 1} & {\bf 1} & 0 & 0 & 1 & 0 & 0 & 1 & 0 & 0 & 1 & 0 & 0 & 1 & 0 & 0 & 1 & 0 & 0 \\ \hhline{-----||-----||------------}
\end{tabular}
\end{table}



\section{Conclusions}

We demonstrated that clausal proof minimization can be an effective technique to reduce the size of graphs with a given property. 
We used this method to shrink graphs while preserving the chromatic number. This resulted in a dozen of
unit-distance graphs with chromatic number 5 consisting of 553 vertices --- a reduction of over
1000 vertices compared to the smallest previously known unit-distance graph with chromatic number 5.

A main goal of this research is to obtain a human-understandable unit-distance graph with 
chromatic number 5. Although that goal has not been reached yet, the experiments produced some interesting results.
For example, either all vertices at distance $\frac{\sqrt{33} - 1}{2\sqrt{3}}$
or all vertices at distance $\frac{\sqrt{33} + 1}{2\sqrt{3}}$ from the central vertex are forced to the same color as the central vertex by the
large part of the minimized graphs.
Also, our research produced a symmetric graph of 199 vertices that was vital for the reduction. We will study
this graph in more detail to determine which properties make it so useful. Moreover, we found two rotations that connected
points at multiple distances, thus increasing the average vertex degree of unit-distance graphs. Finding
more such rotations may allow us to shrink the graphs even further.

Applying clausal-proof techniques to provide mathematical insights is an interesting twist in the discussion about
the usefulness of mechanized mathematics. It has been argued that computers are just ``ticking off possibilities''~\cite{Lamb}.
In this case, however, they reveal important patterns. The techniques described in this paper may actually produce
the most clean and compact proof that the chromatic number of the plane is at least 5. 

Finally, all graphs used in our experiments could be easily colored with 5 colors, even the ones with many thousands of vertices.
However, we observed that this does not hold for the graph $(S_{199} \oplus S_{199}) \cup \theta_4(S_{199} \oplus S_{199})$.
This graph is 5-colorable, even when requiring two colors for the central vertex, but computing such a coloring is expensive. Consequently,
such colorings may be rare and thus may contain certain patterns. This could point to the existence of unit-distance graphs with chromatic number 6 with thousands of vertices.

\section*{Acknowledgements}

The author thanks Aubrey de Grey, Jasmin Blanchette, Benjamin Kiesl, Victor Marek, and the reviewer for their valuable input on an earlier draft of this paper. 
The author acknowledges the Texas Advanced Computing Center (TACC) at The University of Texas at Austin for providing HPC
resources that have contributed to the research results reported within this paper. 

\bibliography{paper}
\bibliographystyle{splncs03}

\end{document}